\documentclass[11pt,a4paper]{article}
\usepackage[top=3cm, bottom=3cm, left=3cm, right=3cm]{geometry}
\usepackage[latin1]{inputenc}
\usepackage{amsmath}
\usepackage{amsthm}
\usepackage{amsfonts}
\usepackage{amssymb}
\usepackage{graphicx}

\usepackage{amsmath,amsthm,amssymb,graphicx, multicol, array}
\usepackage{enumerate}
\usepackage{enumitem}
\usepackage{hyperref}

\usepackage{tabto}

\DeclareMathOperator{\li}{li}

\newtheorem{thm}{Theorem}[section]
\newtheorem{lem}{Lemma}[section]
\newtheorem{conj}{Conjecture}[section]

\newtheorem{dfn}{Definition}[section]

\newcommand{\N}{\mathbb{N}}

\newcommand{\Q}{\mathbb{Q}}

\usepackage{fancyhdr}
\title{Euler Totient Function Inequality}
\date{}
\author{N. A. Carella}

\begin{document}
\thispagestyle{empty}
\maketitle

\begin{abstract}
A new unconditional inequality for the Euler totient function is contributed to the literature. This result is associated with various problems about the distribution of prime numbers. \let\thefootnote\relax\footnote{ \today \date{} 
	}
\end{abstract}

\tableofcontents

\section{Introduction} \label{s1}\hyperget{s1}
The Euler totient function $\varphi(n)=\#\{ m < n : \gcd(m,n) = 1 \}$, which counts the number of relatively prime integers less than $n$, is a sine qua non in number theory. It and its various generalizations appear everywhere in the mathematical literature. The product form representation
\begin{equation}
\varphi(n)=\prod_{p \mid n} \left( 1-\frac{1}{p} \right )
\end{equation}
unearths its intrinsic link to the distribution of the prime numbers. The totient function $\varphi(n)$ is an oscillatory function, its value oscillates from its maximum $\varphi(n)=n-1$ at prime integers $n$ to its minimum $\varphi(n)=c_0 n/ \log \log n$, at the primorial integers $N_k= 2\cdot 3 \cdots p_k$, where $p_i$ is the $i$th
prime in increasing order, and $c_0 > 0$ is a constant. The new contributions to the literature are the unconditional estimates stated below.

\begin{thm} \label{thm1.2A}\hyperget{thm1.2A} Let $n \geq1$ be a large integer, then 
\begin{enumerate} [font=\normalfont, label=(\roman*)]
\item $\displaystyle \frac{n}{\varphi(n)}\leq 
e^{\gamma} \log \log n + O\left(\frac{(\log \log n)^2}{(\log n)^{1/8}} \right),$\tabto{9cm}Unconditionally.\\

\item$\displaystyle \frac{\varphi(n)}{n}\geq 
\frac{1}{e^{\gamma} \log \log n} + O\left(\frac{\log \log n}{(\log n)^{1/8}} \right),$\tabto{9cm}Unconditionally.
\end{enumerate}
\end{thm}
The proof is based on a new unconditional zerofree region of the zeta function, called the quasi-RH, appears in \hyperlink{s7}{Section} \ref{s7}. The above improves the current best unconditional estimate of this arithmetical function in the literature. 
\begin{thm} \label{thm1.2}\hyperget{thm1.2} {\normalfont (\cite{AC2024}) } Let $n \in \N$, then 
\begin{equation}
\frac{n}{\varphi(n)} \leq 
e^{\gamma} \log \log n +O\left(e^{-c(\log \log k)^{3/5}}\right)
	\end{equation}
where $a>0$ is a constant, unconditionally.
\end{thm} 
The error term is reduced by the subexponential factor
\begin{equation}
\frac{(\log\log n)^2e^{-(\log\log n)/8}}{e^{-a(\log\log k)^{3/5}}}=(\log\log n)^2e^{-(\log\log n)/8+a(\log\log k)^{3/5}}.
\end{equation}

The lower and upper bounds for almost every integer are significantly smaller by an iterated factor of log as demonstrated below.
\begin{thm} \label{thm1.4}\hyperget{thm1.4} For almost all integers $n \geq 1$, the ratio $n/\varphi(n)$ has the followings bounds.
	\begin{enumerate} [font=\normalfont, label=(\roman*)]
		\item $\displaystyle  \log \log \log n \ll \frac{n}{\varphi(n)},$\tabto{9cm}Unconditionally.\\
		
		\item$\displaystyle \frac{n}{\varphi(n)} \ll \log \log \log n,$\tabto{9cm}Unconditionally.
	\end{enumerate}
\end{thm}

\begin{thm} \label{thm1.5}\hyperget{thm1.5}  Let $ x \geq 1$ be a large number. Then, for almost all large integers $n \geq 1$, the Euler totient function has the followings bounds.
	\begin{enumerate} [font=\normalfont, label=(\roman*)]
		\item $\displaystyle \frac{\varphi(n)} {n} \geq \frac{c_0}{\log \log \log n} 
		,$\tabto{6cm}for some constants $c_0>0$, unconditionally.\\
		
		\item$\displaystyle \frac{\varphi(n)} {n} \leq \frac{c_1}{\log \log \log n},$\tabto{6cm}for some constants $c_1>0$, unconditionally.
	\end{enumerate}
\end{thm}
On the other hand, there are several conditional criteria; one of these is listed below.

\begin{thm} \label{thm1.3}\hyperget{thm1.3} {\normalfont (\cite{NJ83}) } Let $n_k= 2\cdot 3\cdots p_k$ be the product of the first $k\geq 1$ primes.
\begin{enumerate} [font=\normalfont, label=(\roman*)]
\item If the Riemann Hypothesis is true, then, for each $n_k \geq 5041,$
$$
\frac{n_k}{\varphi(n_k)} > e^{\gamma} \log \log n_k
$$ for all $k \geq 1$.
\item If the Riemann Hypothesis is false, then, 
$$
\frac{n_k}{\varphi(n_k)} <e^{\gamma} \log \log n_k \quad \text{ and } \quad \frac{n_k}{\varphi(n_k)} > 
e^{\gamma} \log \log n_k$$
occur for infinitely many $k \geq 1$.
\end{enumerate}
\end{thm}
\begin{thm} Under the Riemann hypothesis (RH), the followings hold.
	\begin{enumerate} [font=\normalfont, label=(\roman*)]
		\item $\displaystyle \limsup_{n\to\infty}	c(n) = e\gamma(2 +\beta) = 3.6444\ldots,$
		\item $\displaystyle c(n) < e\gamma(2 +\beta),$\tabto{7.7cm}for all $n\geq N_{120569} = 2\cdot3\cdots1591883$,
		\item $\displaystyle c(n) \leq c(N_{66}) = c(2\cdot 3\cdots 317) = 4.0628\ldots,$\tabto{7.7cm}for all $n\geq N_{120569} = 2\cdot3\cdots1591883$,
		\item $\displaystyle c(N_k) \geq c(N_1) = c(2) =2.2086\ldots,$\tabto{7.7cm}for all $n\geq 1$.
	\end{enumerate}	
\end{thm}
\begin{proof}The proof appears in \cite[Theorem 1.1]{NS12}.
\end{proof}

Some related and earlier works on this topic include the works of Ramanujan, Erdos, and other on abundant numbers, see \cite{RS15}, \cite{AE44}, and recent related works appeared in \cite{BK06}, \cite{BN07}, \cite{LJ00}, \cite{NS12}, and \cite{WM07}. The first few \hyperlink{s200}{Section} \ref{s200} to \hyperlink{s5}{Section} \ref{s5} cover some background materials focusing on some finite sums over the prime numbers and some associated and products. The proofs of \hyperlink{thm1.3}{Theorem} \ref{thm1.3}, \hyperlink{thm1.4}{Theorem} \ref{thm1.4} and \hyperlink{thm1.5}{Theorem} \ref{thm1.5} respectively are given in the last \hyperlink{s8}{Section} \ref{s8}.

\section{Prime Numbers Theorems} \label{s200}\hyperget{s200}
The omega notation $f (x)= g(x)+\Omega_{\pm}(h(x))$ means that both $f (x) > g(x)+ c_0h(x)$ and $f (x) < g(x)- c_1h(x)$ occur infinitely often as $x \to \infty$, where $c_0 > 0$ and $c_1 > 0$ are constants, see \cite[p.\ 5]{MV07}, and similar references.\\

The weighted primes counting functions, psi $\psi(x)$ and theta $\theta(x)$, are defined by
\begin{equation}
\theta(x)=\sum_{p \leq x} \log p \quad \text{ and } \quad \psi(x)=\sum_{p^k \leq x} \log p^k                 
\end{equation}
respectively. The standard prime counting function is denoted by
\begin{equation}
\pi(x)=\#\{p \leq x\}=\sum_{p \leq x} 1.                 
\end{equation}
This function is usually expressed in term of the logarithm integral $\li(x)=\int_2^x (\log t)^{-1} dt$.

\begin{thm} \label{thm200.1}\hypertarget{thm200.1} Uniformly for $x \geq 2$ the psi and theta functions have the followings asymptotic formulae.
\begin{enumerate} [font=\normalfont, label=(\roman*)]
\item  \text{Unconditionally,}
$$
\theta(x)=x +O\left (x^{7/8}\log x\right ) .
$$
\item \text{Unconditional oscillation,}
$$
\theta(x)=x +\Omega_{\pm} \left ( x^{1/2}\log \log \log x \right ).  $$

\item \text{Conditional on the RH,}
$$
\theta(x)=x +O\left (x^{1/2} \log^2 x \right ).  
$$
\end{enumerate}
\end{thm}

\begin{proof} (ii) The oscillations form of the theta function is proved in \cite[p.\ 479]{MV07}, 
\end{proof}

The same asymptotics hold for the function $\psi(x)$. Explicit estimates for both of these functions are given in \cite{CP85}, \cite{SL76}, \cite[Theorem 5.2]{DP10}, and related literature. 

\begin{conj} \label{conj200.1} Assuming the RH and the LI conjecture, the suprema are
\begin{equation}
\lim \inf_{x \to \infty} \frac{\psi(x)-x}{\sqrt{x} (\log \log x)^2}=\frac{-1}{\pi} \qquad \text{and} \qquad \lim \sup_{x \to \infty} \frac{\psi(x)-x}{\sqrt{x} (\log \log x)^2}=\frac{1}{\pi}.
\end{equation}
\end{conj}

More details on the Linear Independence conjecture appear in \cite{IA42}, \cite[Theorem 6.4]{EE85}, and recent literature. The LI conjecture asserts that the imaginary parts of the nontrivial zeros $\rho_n=1/2+i \gamma_n$ of the zeta function $\zeta(s)$ are linearly independent over the set $\{-1,0, 1\}$. In short, the equations  
\begin{equation}
\sum_{1 \leq n \leq M} r_n \gamma_n=0,
\end{equation}
where $r_n \in \{-1,0, 1\}$, have no nontrivial solutions. 

\begin{thm} \label{thm200.2} Let $x \geq 1$ be a large number. Then
\begin{enumerate} [font=\normalfont, label=(\roman*)]
\item  \text{Unconditionally,}
$$
\pi(x)=\frac{x}{\log x} +O\left (x^{7/8}\log x\right ) .
$$
\item \text{Unconditional oscillation,}
$$
\pi(x)=\li(x) +\Omega_{\pm} \left (\frac{ x^{1/2}\log \log \log x}{\log x} \right ) . $$

\item \text{Conditional on the RH,}
$$
\pi(x)=\li(x) +O\left (x^{1/2} \log x \right ).  
$$
\end{enumerate}
\end{thm}

\begin{proof} (i) The unconditional part of the prime counting formula arises from the delaVallee Poussin form $\pi(x)=\li(x)+O\left (xe^{-c_0 \sqrt{\log x}}\right )$ of the prime number  theorem, see \cite[p.\ 179]{MV07}. Recent information on the constant $c_0>0$ and the sharper estimate $\pi(x)=\li(x)+O\left ( x e^{-c_0 \log x^{3/5}(\log \log x)^{-2/5}}\right )$ appears in \cite{FK00}. The new result follows from the work in \cite[Theorem 1]{OI2026}.\\

(ii) The unconditional oscillations part arises 
from the Littlewood form $\pi(x)=\li(x)+\Omega_{\pm} \left (x^{1/2}\log \log \log  x /\log x \right )$ of the prime number theorem, 
consult \cite[p.\ 51]{IA03}, \cite[p.\ 479]{MV07}, et cetera. \\

(iii)  The conditional part arises from the Riemann form $\pi(x)=\li(x)+O\left (x^{1/2}\log^2 x \right )$ of the prime number theorem. 
\end{proof}

New explicit estimates for the number of primes in arithmetic progressions are computed in \cite{BR18}.
\section{Harmonic Sums over the Primes} \label{s2} \hyperget{s2}
The most basic finite sum over the prime numbers is the prime harmonic sum $\sum_{p \leq x}1/p$. The refined estimate of this finite sum, stated below, is a synthesis of various results due to various authors. 

\begin{lem} \label{lem2.1}\hyperget{lem2.1}
Let \(x\geq 2\) be a large number, then
\begin{enumerate} [font=\normalfont, label=(\roman*)]
\item  $\displaystyle\sum_{p \leq x} \frac{1}{p}
=\log \log x+ B_1 +O\left (\frac{\log x}{x^{1/8}} \right ), $ \tabto{8cm}  Unconditionally,
\item $\displaystyle \sum_{p \leq x} \frac{1}{p}
=\log \log x+ B_1 +\Omega_{\pm} \left (\frac{\log \log \log x}{x^{1/2}\log x} \right ),  $\tabto{8cm} Unconditional oscillation
\item $\displaystyle \sum_{p \leq x} \frac{1}{p}
=\log \log x+ B_1 +O\left (\frac{\log x}{ x^{1/2}} \right ),$\tabto{8cm}Conditional on the RH,
\end{enumerate}
where $B_1 = 0.2614972128 \ldots, $ is Mertens constant, and $c_0>0$ is an absolute constant. 
\end{lem}

\begin{proof} (i) Replace the asymptotic formula $\pi(x)=\li(x)+O(x^{7/8}\log x)$, and the appropriate prime counting measure $\pi(x)$ in Theorem \ref{thm200.2} into the Stieltjes integral representation and evaluate it.
\begin{align}
\sum_{p \leq x}  \frac{1}{p}&=\int_2^{x}\frac{1}{t}d\pi(t) \\[.3cm]
&=\frac{\li(x)+O(x^{7/8}\log x)}{x}+\int_2^{x}\frac{\li(t)+O(t^{7/8}\log t)}{t^2}dt \\[.3cm]
&=\log\log x+B_1+O\left(\frac{\log x}{x^{1/8}} \right), 
\end{align}
where $B_1>0$ is a constant.\\

(ii) The unconditional oscillations part arises 
from the Littlewood form $\pi(x)=\li(x)+\Omega_{\pm} \left (x^{1/2}\log \log \log  x /\log x \right )$ of the prime number theorem, consult \cite[p.\ 51]{IA03}, \cite[p.\ 479]{MV07}, et cetera. \\

(iii)  The conditional part arises from the Riemann form $\pi(x)=\li(x)+O\left (x^{1/2}\log^2 x \right )$ of the prime number theorem. 
\end{proof}

The asymptotic order $\sum_{p \leq x}1/p \sim \log \log x$ is due to Euler, confer \cite[Chapter 15?]{EL88}. The earliest version including error term $\sum_{p \leq x}1/p=\log \log x+ B_1 +O(1/ \log x)$ is due to Mertens, see \cite{VM05}. The qualitative form of the oscillations of the differences 
\begin{equation} \label{eq200.8}
\sum_{p^k \leq x} \frac{1}{p^k}-\left ( \log \log x + \gamma \right ) \quad \text{ and} \quad \sum_{p^k \leq x} \frac{\log p}{p^k}-\left ( \log x + \gamma \right )
\end{equation} 
seems to be due to Phragmen, confer \cite[p.\ 182]{NW00}. \\

The Euler constant and Mertens constant occur very frequently in analysis. The former is defined by 
\begin{equation}
\gamma=\lim_{x \to \infty} \left ( \sum_{n \leq x}\frac{1}{n}-\log x \right ) =0.577215665 \ldots,
\end{equation}
and the later 
is defined by 
\begin{equation}
B_1=\lim_{x \to \infty} \left ( \sum_{p \leq x}\frac{1}{p-1}-\log \log x \right ) =0.2614972128 \ldots.
\end{equation}
Other definitions of these constants are available in the literature, confer \cite{LJ13}.
 
\begin{lem} \label{lem2.2} The constants $\gamma$ and $B_1$ satisfy the linear relation
\begin{equation}
B_1=\gamma-\sum_{p \geq 2} \sum_{n \geq 2} \frac{1}{np^{n}}.                 
\end{equation}

\end{lem}

\begin{proof} This relation stems from the power series expansion 
\begin{equation}
B_1-\gamma=\sum_{p \geq 2} \left (\log \left (1-\frac{1}{p} \right ) + \frac{1}{p} \right )
\end{equation}
via the power series for $\log(1+z)$ with $|z|<1$. This leads to this identity, see \cite[p.\ 466]{HW79}, \cite[p.\ 182]{MV07}.\end{proof}

\begin{lem} \label{lem2.3} Let \(x\geq 2\) be a large number, then
\begin{equation}
\sum_{p \leq x} \sum_{n \geq 2} \frac{1}{np^{n}} =\gamma-B_1-\frac{1}{x \log x}+ O \left (\frac{1}{x\log^2 x}\right ). 
\end{equation}
\end{lem}

\begin{proof} Rearrange the power series expansion as
\begin{equation}
\sum_{p \leq x} \sum_{n \geq 2} \frac{1}{np^{n}} =\gamma-B_1- \sum_{p > x} \sum_{n \geq 2} \frac{1}{np^{n}} =\gamma-B_1-\frac{1}{x \log x}+ O \left (\frac{1}{x\log^2 x}\right ). 
\end{equation}
The estimate for the last two terms on the right follows from Lemma \ref{lem2.4} computed below.
\end{proof} 

\begin{lem} \label{lem2.4}  Let $x \geq 1$ be a large number, then
\begin{equation}
\sum_{p >x}\sum_{n \geq 2} 
\frac{1}{np^{n}} = \frac{1}{x \log x}+ O \left (\frac{1}{x\log^2 x}\right ).
\end{equation}
\end{lem}
\begin{proof} Split the infinite sum into two subsums:
\begin{eqnarray}
\sum_{p > x} \sum_{n \geq 2} \frac{1}{np^n}&=&\sum_{p > x}  \frac{1}{2p^2}+\sum_{p > x} \sum_{n \geq 3} \frac{1}{np^n} \nonumber \\
&=&\sum_{p \geq x}  \frac{1}{2p^2}+O \left (\frac{1}{x^2 \log x}\right) \nonumber.
\end{eqnarray}
Employ the prime counting measure $\pi(t)=\#\{p \leq t\}$ to evaluate the first subsum using the integral 
\begin{eqnarray}
\sum_{p \geq x}  \frac{1}{p^2}&=&\int_x^{\infty}\frac{1}{t^2}d\pi(t) \nonumber \\
&=&-\frac{\pi(x)}{x^2} +2\int_x^{\infty}\frac{\pi(t)}{t^3}dt \nonumber \\
&=&\frac{1}{x \log x}+O \left (\frac{1}{x \log^2 x}\right) \nonumber .
\end{eqnarray} 
\end{proof}

\section{Products over the Primes} \label{s3} \hyperget{s3}
The asymptotics for a variety of interesting products are simple applications of the results for prime harmonic sums in the previous section.
\begin{lem} \label{lem300.1}\hypertarget{lem300.1}
Let \(x\geq 2\) be a large number, then

\begin{enumerate} [font=\normalfont, label=(\roman*)]
\item  $\displaystyle \prod_{p \leq x}\left( 1- \frac{1}{p} \right) ^{-1}
=e^{\gamma} \log x+ O\left (\frac{\log x}{x^{1/8}} \right ) ,$\tabto{9cm}Unconditionally,
\item $\displaystyle \prod_{p \leq x}\left( 1- \frac{1}{p} \right) ^{-1}
=e^{\gamma} \log x+\Omega \left (\frac{\log \log \log x}{x^{1/2}} \right ),$\tabto{9cm}Unconditional oscillation,
\item $\displaystyle \prod_{p \leq x}\left( 1- \frac{1}{p} \right) ^{-1}
=e^{\gamma} \log x+ O\left (\frac{\log x}{ x^{1/2}} \right ) ,$\tabto{9cm} Conditional on the RH,
\end{enumerate}
where $B_1 = 0.2614972128 \ldots, $ is Mertens constant, and $c_0>0$ is an absolute constant. 
\end{lem}
\begin{proof}(i) This follows from \hyperlink{lem2.1}{Lemma}  \ref{lem2.1}.
\end{proof}

The results for products over arithmetic progression are proved in \cite{LZ07}, et alii.\\

The nonquantitative unconditional oscillations of the error of the product of primes is implied by the work of Phragmen, refer to equation (\ref{eq200.8}), and \cite[p.\ 182]{NW00}. Since then, various authors have developed quantitative versions, see \cite{RS62}, \cite{DP09}, \cite{LY14}, \cite{LJ15}, et alii. The specific quantitative form
\begin{equation} \label{eq3.10}
\prod_{p \leq x}\left( 1- \frac{1}{p} \right) ^{-1}
=e^{\gamma} \log x+ \Omega_{\pm}\left (\frac{f(x)}{ x^{1/2}} \right ), 
\end{equation}
where $f(x)$ is a slowly increasing function, was proved in \cite{DP09}. \\

\begin{thm} \label{thm200.21} { \normalfont (Martens, 1874)} The following asymptotic formulas hold:
\begin{enumerate} [font=\normalfont, label=(\roman*)]
\item  
$$
\lim_{x \to \infty} \frac{1}{\log x} \prod_{p \leq x} \left (1-\frac{1}{p}\right )^{-1}=e^{\gamma} .
$$
\item 
$$
\lim_{x \to \infty} \frac{1}{\log x} \prod_{p \leq x} \left (1+\frac{1}{p}\right )^{-1}=\frac{6e^{\gamma}}{\pi^2} .
$$
\end{enumerate}
\end{thm}

\section{Logarithm Sums of Primes} \label{s4}\hyperget{s4}
The logarithm of a product can be derived from the estimate of the product itself. However, an independent proof based on the power series of the logarithm will be used here to obtain these estimates.
\begin{lem} \label{lem600.1}  Let $x \geq 1$ be a large number, then
\begin{enumerate} [font=\normalfont, label=(\roman*)]
\item 
$\displaystyle \log \prod_{p \leq x}\left( 1- \frac{1}{p} \right) ^{-1}
=\log \log x+\gamma  + O\left (\frac{\log x}{x^{1/8}} \right ),$\tabto{9.5cm} Unconditionally,
\item $\displaystyle \log \prod_{p \leq x}\left( 1- \frac{1}{p} \right) ^{-1}
=\log \log x+\gamma  +\Omega \left (\frac{\log \log \log x}{x^{1/2} \log x} \right ),$\tabto{9.5cm}Unconditional oscillation,
\item$\displaystyle \log \prod_{p \leq x}\left( 1- \frac{1}{p} \right) ^{-1}
=\log \log x+\gamma  + O\left (\frac{\log \log \log \log x}{x^{1/2}\log x} \right ),$\tabto{10cm}Conditional on the RH,
\end{enumerate}
where $p \leq x$ is the largest prime divisor of $n$.
\end{lem}

\begin{proof} (ii) The product form of the ratio $n/\varphi(n)$ has the equivalent form
\begin{eqnarray}
 \prod_{p \leq x }\left( 1- \frac{1}{p} \right) ^{-1}  &=&  \prod_{p \leq x}\left( 1- \frac{1}{p} \right) ^{-1} \left( 1+ \frac{1}{p} \right) ^{-1}\left( 1+ \frac{1}{p} \right)\nonumber \\
&=& \prod_{p \leq x}\left( 1- \frac{1}{p^2} \right) ^{-1} \left( 1+ \frac{1}{p} \right) .
\end{eqnarray}
Taking logarithms and replacing the power series expansions return
\begin{eqnarray}
\sum_{p \leq x} \left( \log \left( 1-\frac{1}{p^2} \right )^{-1} +\log \left( 1+\frac{1}{p} \right ) \right )&=&\sum_{p \leq x}\sum_{n \geq 1} \frac{1}{np^{2n}} +  \sum_{p \leq x} \sum_{n \geq 1} \frac{(-1)^{n+1}}{np^n}  \\
&=&\sum_{p \leq x}\sum_{n \geq 1} \frac{1}{np^{2n}} +  \sum_{p \leq x} \sum_{n \geq 2} \frac{(-1)^{n+1}}{np^n} +\sum_{p \leq x}  \frac{1}{p}\nonumber \\
&=&\sum_{p \leq x} 
\left (\sum_{n \geq 1} \frac{1}{np^{2n}} +   \sum_{n \geq 2} \frac{(-1)^{n+1}}{np^n}\right ) +\sum_{p \leq x}  \frac{1}{p} \nonumber.
\end{eqnarray}
By\hyperlink{lem2.1}{Lemma} \ref{lem2.1}, the last finite sum
\begin{eqnarray}
\sum_{p \leq x} \frac{1}{p}&=&\log \log x+B+\Omega_{\pm} \left (\frac{\log \log \log x}{ x^{1/2}\log x}\right)  \\
&=&\log \log x+\gamma-\sum_{p \geq 2}\sum_{n \geq 2} \frac{1}{np^{n}}  +\Omega_{\pm} \left (\frac{\log \log \log x}{ x^{1/2}\log x}\right)\nonumber.
\end{eqnarray}
The last line follows from \hyperlink{lem2.3}{Lemma} \ref{lem2.3}. Replacing back and combining the infinite series, yield
\begin{eqnarray}
&=&\sum_{p \leq x} 
\left (\sum_{n \geq 1} \frac{1}{np^{2n}} +   \sum_{n \geq 2} \frac{(-1)^{n+1}}{np^n}\right ) +\sum_{p \leq x}  \frac{1}{p}  \\
&=& \sum_{p \leq x} 
\left (\sum_{n \geq 1} \frac{1}{np^{2n}} +   \sum_{n \geq 2} \frac{(-1)^{n+1}}{np^n}\right ) + \log \log x+\gamma-\sum_{p \geq 2}\sum_{n \geq 2} \frac{1}{np^{n}}  +\Omega_{\pm} \left (\frac{\log \log \log x}{ x^{1/2}\log x}\right)\nonumber\\
&=& \sum_{p \leq x} \left (\sum_{n \geq 1} \frac{1}{np^{2n}} +   \sum_{n \geq 2} \frac{(-1)^{n+1}}{np^n} -\sum_{n \geq 2} \frac{1}{np^{n}}\right ) \nonumber \\
&& \qquad \qquad   \qquad  \qquad + \log \log x+\gamma-\sum_{p >x}\sum_{n \geq 2} \frac{1}{np^{n}}  +\Omega_{\pm} \left (\frac{\log \log \log x}{ x^{1/2}\log x}\right) \nonumber \\
&=& \log \log x+\gamma-\sum_{p >x}\sum_{n \geq 2} \frac{1}{np^{n}}  +\Omega_{\pm} \left (\frac{\log \log \log x}{ x^{1/2}\log x}\right) \nonumber.
\end{eqnarray}
The triple series on the third line vanish, and the series on the fourth line is small enought to be absorved into the omega term. 
 \end{proof}

\begin{lem} \label{lem600.2}  Let $n \geq 1$ be a large highly composite integer, then
\begin{equation}
\log \left ( \frac{n}{\varphi(n)} \right ) = \log \log \log n +\gamma+\Omega_{\pm} \left (\frac{\log \log \log \log n}{ (\log n)^{1/2}\log \log n}\right).
\end{equation}
\end{lem}

\begin{proof} Given a highly composite integer, the magnitude of the largest prime divisor $p \mid n$ is determined in Lemma \ref{lem500.1}, which gives
\begin{equation}
p \leq x=\log n \left (1+\Omega_{\pm}\left (\frac{\log \log \log \log n}{(\log n)^{1/2}} \right )  \right ), 
\end{equation}
and  
\begin{eqnarray} \label{eq999} 
\log \log x&=&\log \log \left (\log n \left (1+\Omega_{\pm}\left (\frac{\log \log \log \log n}{(\log n)^{1/2}} \right ) \right )\right ) \nonumber \\
&=&\log \left (\log  \log n +\log \left (1+\Omega_{\pm}\left (\frac{\log \log \log \log n}{(\log n)^{1/2}} \right ) \right )\right ) \\
&=&\log \log \log n +\Omega_{\pm}\left (\frac{\log \log \log \log n}{(\log n)^{1/2} \log \log n} \right ) \nonumber.
\end{eqnarray}
By the prime number theorem, the error term in (\ref{eq999}) changes sign infinitely often. \\

Applying Lemma \ref{lem600.1}, and replacing (\ref{eq999}), yield
\begin{eqnarray}
\log \left ( \frac{n}{\varphi(n)} \right ) &=&\log \log x+\gamma  +\Omega_{\pm} \left (\frac{\log \log \log x}{ x^{1/2}\log x}\right)  \\
&=& \log \log \log n +\gamma+\Omega_{\pm} \left (\frac{\log \log \log \log n}{ (\log n)^{1/2} \log \log n}\right) \nonumber,
\end{eqnarray}
this proves the claim.
 \end{proof}


\section{Highly Composite Numbers}  \label{s5}\hyperget{s5}
Let $p_k$ be the $k$th prime in increasing order, and let $v_p=\max \{ m: p^m \mid n\}$ is the $p$-adic valuation. Extremely abundant integers, and colossally abundant integers are related to the primorial integers $n=2^{v_2} \cdot 3^{v_3} \cdots p^{v_p}$, but the exponents have certain multiplicative structure $1 \leq v_p\leq \cdots \leq v_3 \leq v_2$. 
\begin{dfn} \label{dfn500.4} {\normalfont Let $d(n)=\sum_{d \mid n}1$. An integer $n \in \N$ is called highly composite if and only if $ d(m) < d(n)$ for all integers $m <n$.}
\end{dfn}

\begin{dfn} \label{dfn500.5} {\normalfont Let $\sigma(n)=\sum_{d \mid n}d$. An integer $n \in \N$ is called colossally abundant if and only if
\begin{equation} \frac{\sigma(m)}{m^{1+\varepsilon} }< \frac{\sigma(n)}{n^{1+\varepsilon}}
\end{equation}
for all integers $m <n$, and some small number $\varepsilon>0$}
\end{dfn}
These numbers are studied in \cite{RS15}, \cite{AE44}, \cite{LJ00}, \cite{BK06}, \cite{NS12}, et alii. \\

\begin{lem} \label{lem500.1} {\normalfont (\cite[Theorem 2]{AE44})} Let $n \geq $ be a large highly composite integer, then
\begin{enumerate} [font=\normalfont, label=(\roman*)]
\item Unconditionally, the largest prime divisor $p \mid n$ has the asymptotic 
$$ p= (\log n)\left (1+O\left (\frac{1}{(\log \log n)^2} \right ) \right ).$$
\item  Unconditionally, the largest prime divisor $p \mid n$ has the asymptotic oscillations
$$ p= (\log n)\left (1+\Omega_{\pm}\left (\frac{\log \log \log n}{(\log n)^{1/2}} \right ) \right ).$$
\item  Modulo the Riemann hypothesis, the largest prime divisor $p \mid n$ has the asymptotic 
$$ p= (\log n)\left (1+O\left (\frac{\log \log n}{(\log n)^{1/2}} \right ) \right ).$$
\end{enumerate}
\end{lem}

\begin{proof} The asymptotic part $p \sim \log n$ follows from Theorem 2 in \cite{AE44}, and the three forms of the error term $R(n)$ follows from 

Theorem \ref{thm200.1}. 
\end{proof}

\section{The Main Result} \label{s7}\hyperget{s7}
Let $N_k=2\cdot3\cdots p_k$ be a primorial number, it is the product of the first $k$ primes $p_k$ in increasing order. The theta function has the asymptotic formula $\theta(x)=\sum_{p\leq x}\log p=x+o(x)$. At $x=N_k$, it has the form
\begin{equation}\label{eq7.000}
\theta(p_k)=\sum_{p\leq p_k}\log p=\log N_k.
\end{equation}

\begin{proof} (\hyperlink{thm1.2A}{Theorem} \ref{thm1.2A}) Let $n\leq N_k$ be an arbitrary integer. Next, utilize  \hyperlink{lem300.1}{Lemma} \ref{lem300.1} to evaluate the product representation.
	\begin{align}\label{eq7.200}
		\frac{n}{\varphi(n)}&=\prod_{p\mid n}\left(1-\frac{1}{p} \right)^{-1}\\[.3cm]
		&\leq\prod_{p\mid N_k}\left(1-\frac{1}{p} \right)^{-1}\nonumber\\[.3cm]
		&\leq\prod_{p\ll \log N_k}\left(1-\frac{1}{p} \right)^{-1}\nonumber\\[.3cm]
		&=  e^{\gamma} \log \log N_k+O\left( \frac{\log \log N_k}{( \log N_k)^{1/8}}\right) \nonumber.
	\end{align}
	This completes the verification.
\end{proof}

\section{Average Lower And Upper Bounds} \label{s8}\hyperget{s8}
The subset of integers that satisfy inequality \eqref{eq7.200} is a very thin subset of integers. The counting function has the asymptotic
\begin{equation}
	K(x)=\li(\log x)+O\left(  \frac{\log x}{(\log\log x)^{B+1}}\right),  
\end{equation}
where $B>0$see \cite[p.\;2]{AC2024}. In contrast, almost every integers has a large lower and upper bounds.

\begin{proof} (\hyperlink{thm1.3}{Theorem} \ref{thm1.3}) By the Erdos-Kac theorem, (confer \cite{EK40}, \cite[Theorem 431]{HW79}, \cite{GS06}, et cetera), almost every integer $n \geq 1$ has $\omega(n) \ll \log \log n$ prime divisors. In addition, the interval $[1, (\log \log n)^2]$ contains $\pi((\log \log n)^2) \gg \log \log n$ primes. Thus, 
\begin{eqnarray}
 \frac{n}{\varphi(n)} &=& \prod_{p \mid n }\left( 1- \frac{1}{p} \right) ^{-1}  \\[.2cm]
&\ll&  \prod_{p \ll (\log \log n)^2}\left( 1- \frac{1}{p} \right) ^{-1} \nonumber \\[.3cm]
&\ll& \log \log \log n \nonumber.
\end{eqnarray}
The reverse inequality is similar.
\end{proof}
The proof of Theorem \ref{thm1.5} uses the same technique.

\section{Problems}
1. Use the identity $B_1=\gamma-\sum_{p \geq 2} \sum_{n \geq 2} (np^{n})^{-1}$ to prove or disprove that $B_1$ and $\gamma$ are linearly independent over the rational numbers $\Q$.
\\

2. Evaluate the finite sum $ \sum_{n \leq x}(n \log n)^{-1+\alpha},$ where $\alpha $ is a real number.



\begin{thebibliography}{999}

\bibitem{AC2024} Axler, C. Inequalities involving the primorial counting function. arXiv:2406.04018


\bibitem{AE44} L. Alaoglu and P. Erdos, On highly composite and similar numbers, Trans. Amer: Math. Soc. 56 (1944),
448-469.
\bibitem{BN07} William D. Banks, Derrick N. Hart, Pieter Moree, C. Wesley Nevans, The Nicolas and Robin inequalities
with sums of two squares, arXiv:0710.2424.
\bibitem{BR18} Michael A. Bennett, Greg Martin, Kevin O'Bryant, Andrew Rechnitzer, Explicit 
bounds for primes in arithmetic progressions, arXiv:1802.00085.
\bibitem{CP85} N. Costa Pereira.  Estimates for the Chebyshev function $\theta(x)$ and $\psi(x)$, Math. Comp., 44(169):211-221, 1985.
\bibitem{DP10} P. Dusart. Estimates of some functions over primes without R.H., arXiv:1002.0442.
\bibitem{DP09} H. G. Diamond and J. Pintz, Oscillation of Mertens product formula. J. Theor. Nombres Bordeaux 21 (2009), no. 3, 523-533.
\bibitem{EE85} Ellison, William; Ellison, Fern, Prime numbers. Wiley-Interscience Publication. New York; Hermann, Paris, 1985.
\bibitem{EL88} Euler, Leonhard. Introduction to analysis of the infinite. Book I. Translated from the Latin and with an introduction by John D. Blanton. 

Springer-Verlag, New York, 1988.
\bibitem{EK40} Erdos, P.; Kac, M. The Gaussian law of errors in the theory of additive number theoretic functions. Amer. J. Math. 62, (1940). 738-742.
\bibitem{FK00} Ford, Kevin Zero-free regions for the Riemann zeta function. Number theory for the millennium, II (Urbana, IL, 2000), 25-56, A K Peters, 

Natick, MA, 2002.
\bibitem{FS03} S.R. Finch, Mathematical constants, Encyclopedia of Mathematics and its Applications 94, (Cambridge
University Press, Cambridge, 2003).
\bibitem{GS06}  Andrew Granville, K. Soundararajan, Sieving and the Erdos-Kac theorem, arXiv:math/0606039. 
\bibitem{HW79} G. H. Hardy and E. M. Wright, An Introduction to the Theory of Numbers, 5th ed., Oxford University
Press, Oxford, 1979. 
\bibitem{IA42} A.E. Ingham, On two conjectures in the theory of numbers. Amer. J. Math., 64 (1942), 313-319.
\bibitem{IA03} Ivic, Aleksandar, The Riemann zeta-function. Theory and applications. Wiley, New York; Dover Publications, Inc., Mineola, NY, 2003.
\bibitem{LJ15} J. P.S. Lay, Sign changes in Mertens first and second theorems, arXiv:1505.03589.
\bibitem{LY14} Y. Lamzouri, A bias in Mertens product formula, arXiv:1410.3777v2. 

\bibitem{LZ07} Alessandro, Languasco, Alessandro, Zaccagnini, On the constant in the Mertens product for arithmetic progressions. I., arXiv:0706.2807. 

\bibitem{MH79} H.L. Montgomery, The zeta function and prime numbers, Proceedings of the Queens Number Theory Conference, 1979, Queens Univ., Ont., 1980, 

1-31.
\bibitem{MV07} Montgomery, Hugh L.; Vaughan, Robert C. Multiplicative number theory. I. Classical theory. Cambridge University Press, Cambridge, 2007.
\bibitem{NJ83} J. L. Nicolas, Petites valeurs de la fonction d Euler, J. Number Theory 17 (1983) 375-388.
\bibitem{NJ97} Jean-Louis Nicolas, Guy Robin. Highly Composite Numbers by Srinivasa Ramanujan, The Ramanujan Journal, 1997, Volume 1, Issue 2, pp 119-

153.
\bibitem{NJ2012} Jean-Louis Nicolas, Small values of the Euler function and the Riemann hypothesis, arXiv:1202.0729, or Acta Arithmetica, 155.3, 2012, 

311-321.
\bibitem{NW00}  Narkiewicz, W. The development of prime number theory. From Euclid to Hardy and Littlewood. Springer Monographs in Mathematics. 

Springer-Verlag, Berlin, 2000. 
\bibitem{RS15} S. Ramanujan, Highly composite numbers, Proc. London Math. Soc. 14 (1915), 347-407.
\bibitem{RS62} J.B. Rosser and L. Schoenfeld, Approximate formulas for some functions of prime numbers, Illinois J.
Math. 6 (1962) 64-94.
\bibitem{SL76} Schoenfeld, Lowell. Sharper bounds for the Chebyshev functions $\theta (x)$ and $\psi (x)$. II. Math.
Comp. 30 (1976), no. 134, 337-360.
\bibitem{VM05} Mark B. Villarino, Mertens' Proof of Mertens' Theorem, arXiv:math/0504289.

\bibitem{BK06} Briggs, K. Abundant numbers and the Riemann hypothesis. Experiment. Math. 15 (2006), no. 2, 251-256.
\bibitem{CS06} Y.-J. Choie, N. Lichiardopol, P. Moree, P. Sole, On Robin's criterion for the Riemann Hypothesis,
arXiv:math/0604314.
\bibitem{GT13} Gronwall, T. H. Some asymptotic expressions in the theory of numbers. Trans. Amer. Math. Soc. 14 (1913), no. 1, 113-122.

\bibitem{LJ00} Jeffrey C. Lagarias, An Elementary Problem Equivalent to the Riemann Hypothesis, arXiv:math/0008177.
\bibitem{LJ13} Jeffrey C. Lagarias, Euler's constant: Euler's work and modern developments, arXiv:1303.1856.

\bibitem{OI2026} OpenAI. Proof of the Quasi-Riemann Hypothesis. Preprint 2026, October 6.


\bibitem{RG84} Guy Robin, Grandes valeurs de la fonction somme des diviseurs et hypothese de Riemann. J. Math. Pures Appl. (9) 63 (1984), 187-213.


\bibitem{NS12} S. Nazardonyavi, Superabundant numbers, their subsequences and the Riemann, arxiv.org/pdf/1211.2147.

\bibitem{WM07} Marek Wojtowicz, Robins inequality and the Riemann hypothesis, Proc. Japan Acad. Ser. A Math. Sci.
Volume 83, Number 4 (2007), 47-49.



 




 
\end{thebibliography}
\end{document}